\documentclass[12pt]{article}

\usepackage[utf8]{inputenc}
\usepackage[T1]{fontenc}

\usepackage{amsthm}
\usepackage{amsmath}
\usepackage{amssymb}
\usepackage{mathrsfs}
\usepackage{stmaryrd}
\usepackage{indentfirst}

\usepackage{graphicx}

\usepackage{hyperref}
\usepackage{lmodern}
\usepackage{fullpage}
\usepackage{fancybox}
\usepackage{xcolor}
\usepackage{tikz}

\usepackage[style=alphabetic]{biblatex}
\usepackage{csquotes}

\theoremstyle{plain}
\newtheorem{thm}{Theorem}[section]
\newtheorem{prop}[thm]{Proposition}

\newtheorem{lem}[thm]{Lemma}
\newtheorem{cor}[thm]{Corollary}
\newtheorem{fact}[thm]{Fact}
\newtheorem*{claim}{Claim}
\newtheorem*{question}{Question}

\theoremstyle{definition}
\newtheorem{defn}[thm]{Definition}
\newtheorem*{exemple}{Example}

\theoremstyle{remark}
\newtheorem*{remark}{Remark}
\newtheorem*{remarks}{Remarks}

\newcommand\N{\mathbb N}

\newcommand\Q{\mathbb Q}
\newcommand\U{\mathcal U}
\newcommand\tp{\textrm{tp}}
\newcommand\Aut{\textrm{Aut}}
\def\forkindep{\mathrel{\raise0.2ex\hbox{\ooalign{\hidewidth$\vert$\hidewidth\cr\raise-0.9ex\hbox{$\smile$}}}}}
\newcommand\mes{\mathfrak M}

\begin{document}

\begin{center}

\huge{Smooth measures and the canonical retraction in NIP theories}

\rule{0.8\textwidth}{0.5mm}

\vspace{1cm}

\small{Xavier Pigé}


\vspace{1.5cm}

\textbf{Abstract}

\vspace{.5cm}

	\begin{minipage}{0.9\textwidth}
	\small{We show that the results proved by Simon on the canonical retraction $F_M$ from the space of $M$-invariant types onto the space of types finitely satisfiable in $M$ remain true over measures. We also make another construction of the canonical retraction for measures, mimicking what Simon did for types, and show that it coincides with Simon's canonical retraction for measures. To do so, we make extensive use of smooth measures.}
\end{minipage}

\end{center}

\pagestyle{plain}

\section*{Introduction}

A \textit{(Keisler) measure} (over $A$) is a finitely additive measure on the Boolean algebra of ($A$-)definable sets. It is \textit{smooth} if it has a unique global extension. When interpreting types as special cases of measures, a type over a model gives rise to a smooth measure if and only if it is realized. In NIP theories, Proposition 7.9 in \cite{NIP_guide} states that any Keisler measure can be extended to a smooth measure, which is analogous to saying that any type is realized in an extension. Using this result as a basis, it is possible to establish a correspondance between types and measures. This correspondance allows to translate results (their statements and their proofs) from types to measures. We will here do this in the special case of the canonical retraction, introduced by Simon in \cite{inv_types_NIP}.

Let us explain what is Simon's canonical retraction. Let $M$ be a small model, and $p$ an $M$-invariant type. Let $M\prec^+M'$ be a $(|M|+|T|)^+$-saturated extension, and $(M',M)\prec^+(N',N)$ (where $(M',M)$ denotes the model $M'$ with an additional unary predicate $P(x)$, interpreted by $P(M')=M$). It follows from NIP that $p_{|N}\cup\{P(x)\}$ implies a complete type over $M'$ (in the original language), say $q$. But then $q$ is finitely satisfiable in $M$, since $q\cup\{P(x)\}$ is consistent. Thus $q$ has a unique global $M$-invariant extension $F_M(p)$, which is also finitely satisfiable in $M$. It can be shown that this is well defined, and that $F_M$ is indeed a continuous retraction from the space of $M$-invariant types onto that of types which are finitely satisfiable in $M$. It is also possible to extend $F_M$ to an affine continuous retraction from the space of $M$-invariant measures onto that of measures which are finitely satisfiable in $M$, by setting $F_M(\mu)=F_{M*}(\mu)$ (with $\mu$ finitely satisfiable in $M$ meaning that, for all $\phi(x)\in L(\U)$, then $\mu(\phi(x))>0$ implies that $\phi(x)$ is satisfiable in $M$; here affine means that if $\mu,\nu$ are two $M$-invariant measures and $t\in[0,1]$, then $F_M((1-t)\mu + t\nu) = (1-t)F_M(\mu)+tF_M(\nu)$). For more details, see Section 3 of \cite{inv_types_NIP} (for the case of types) and Section 3 of \cite{ext_def_NIP_gps} (for the measure case).

In this note, we will provide a few tools for the study of smooth measures in a theory where not all formulas are NIP - but those we are interested in are. We apply this to pairs of models. More specifically, we propose an alternative construction of $F_M$ over the space of measures (in Subsection \ref{ssec_alt_cons}), showing that, actually, it is also possible to build this canonical retraction through the same process as the one used for types.

Another important object in this paper will be definable measures. Recall that a measure is said to be $M$-definable if it is $M$-invariant and the maps $f_\varphi^\mu:q\in S_y(M)\mapsto \mu(\phi(x;b))\in[0,1]$ (where $b$ is any realization of $q$) are continuous for all $\phi(x;y)\in L$. Also recall that there is a product over (global) types, defined as soon as the left argument is invariant (see Subsection \ref{ssec_prod}). It is known that definable types always commute with finitely satisfiable types (see Lemma 2.23 in \cite{NIP_guide}).

In \cite{dp_types} and \cite{inv_types_NIP}, Simon shown the following converses under NIP hypothesis:

\begin{fact}\label{main_thm_type}
Let $p$ be an $M$-invariant type.

 \begin{enumerate}
 \item Assume that, for all $q$ finitely satisfiable in $M$, we have $(p\otimes q)_{|M} = (q\otimes p)_{|M}$. Then $p$ is $M$-definable.
 \item Assume that $p\otimes F_M(p) = F_M(p)\otimes p$. Then $p$ is $M$-definable.
 \end{enumerate}
 \end{fact}
 
In NIP, this product over types can be generalized to measures, which yields a product of measures defined as soon as the left argument is invariant. It is true that definable measures commute with finitely satisfiable measures (see Proposition 7.22 in \cite{NIP_guide}). One of the main results proved here (Theorem \ref{thm_F_M_def}) will be the following, a measure equivalent the previous theorem:

\begin{thm}\label{main_thm}
Let $\mu$ be an $M$-invariant measure.

\begin{enumerate}
 \item Assume that, for all $q$ finitely satisfiable in $M$, we have $(\mu\otimes q)_{|M} = (q\otimes \mu)_{|M}$. Then $\mu$ is $M$-definable.
 \item Assume that $\mu\otimes F_M(\mu) = F_M(\mu)\otimes \mu$. Then $\mu$ is $M$-definable.
 \end{enumerate}
\end{thm}

The paper is organised as follows. We start by studying the heir property for measures (which shares a strong connection with definability), in Section \ref{section_heir}. Section \ref{section_smooth} is dedicated to explaining why smooth measures in NIP work analogously as realized types. In Section \ref{section_constr}, we develop tools for using smooth measures in the context of pairs of models, and then propose an alternative construction for $F_M$ on measures. Finally, in Section \ref{section_def}, we show Theorem \ref{main_thm}, first showing an equivalent to heir-coheir duality and then following Simon's proof in \cite{inv_types_NIP}, but translating it to measures.



\section{Preliminaries}

Throughout all the text, $\mathcal U$ will be the monster model, i.e. a big saturated model, and $L$ is the language in which we work. We will start by recalling some results on Keisler measures, and then on the products (both on types and measures). We will often use the notation $M\prec^+ N$, which means that $N$ is a $(|M|+|T|)^+$-saturated elementary extension of $M$.

\subsection{The space of measures}


Let us start by some results on Keisler measures. If $x$ is a variable and $M\prec\U$ a model (or even any subset of $A$), not necessarily small, we set $\mes_x(M)$ to be the space of Keisler measures in variable $x$ over $M$.

First note that Keisler measures over $M$ happen to be in natural bijection with Borel regular ($\sigma$-additive) probability measures on $S_x(M)$ (see Section 7.1 in \cite{NIP_guide}). We only explain how to go from a Borel regular probability measure $\mu$ towards a Keisler measure $\tilde\mu$. Let $X$ be an $M$-definable set, say $X = \varphi(\U)$ for $\varphi(x)\in L(M)$. Then $\varphi(x)$ defines a clopen set $[\varphi(x)]$ in $S_x(M)$. Just set $\tilde\mu(X) :=\mu([\varphi(x)])$. We most often write $\tilde\mu(\phi(x))$ for $\tilde\mu(X)$, when $X$ is defined by $\phi(x)$, and make no difference between $\tilde \mu$ and $\mu$.

We will often use measures in several variables, say $\lambda(x,y)$. In this case, we write $\lambda_{|x}$ for the measure in variable $x$ naturally attached to $\lambda$, i.e. $\lambda_{|x}(\varphi(x)):=\lambda(\varphi(x))$ for all $\varphi(x)\in L(M)$. Also, we will need to look at restrictions of measures to submodels. For instance, if $\lambda$ is a measure over $M$ and $N\prec M$, there is a measure over $N$ naturally defined by restricting $\lambda$; we denote this measure by $\lambda_{|N}$. If we need to restrict simultaneously the base space and the variables, we write for instance $\lambda_{|x,N}$. Also, we sometimes write $\mu_x$ instead of $\mu$ to say that $x$ is the variable of $\mu$.

We sometimes speak of \textit{partial} Keisler measures, by analogy with partial types. By this we mean an additive probability measure on a Boolean subalgebra of the Boolean algebra of definable subsets of $M$.

Types can be seen as a special case of measure: they are in natural correspondance with $0-1$-valued measures. Sometimes (espacially when manipulating products), we will write types in some places where a measure would be required; this means we consider the measure attached to this type.

\begin{fact}\label{compact_measures}
\begin{itemize}
\item The space $\mes_x(\U)$ of global Keisler measures is compact (as a subset of $[0,1]^{L_x(\U)}$ endowed with the product topology, where $L_x(\U)$ is the set of definable subsets of $\U$).
\item Any partial measure can be extended to a global Keisler measure.
\item Moreover, if $\mu$ is a partial measure and $\phi(x)\in L(\U),r\in[0,1]$, then $\mu$ has a global extension satisfying $\mu(\phi(x))=r$ if and only if for all $\psi(x)$ measured by $\mu$, we have neither $\psi(x)\to\phi(x)$ and $\mu(\psi(x))>r$, nor $\phi(x)\to\psi(x)$ and $\mu(\psi(x))<r$. 
\end{itemize}
\end{fact}
The first point is rather obvious; the second and third are respectively Lemmas 7.3 and 7.4 in \cite{NIP_guide}.

%

Finally, let us also recall some facts from \cite{gen_stab_smooth} about indiscernible measures:

\begin{defn}
A measure $\mu((x_i)_{i\in I})$, where $I$ is an infinite linear order, is \textit{indiscernible over $M$} if for any $\phi(x_0,\dots,x_{n-1})\in L(M)$ and any $i_0<\dots<i_{n-1},j_0<\dots<j_{n-1}\in I$, we have $\mu(\phi(x_{i_0},\dots,x_{i_{n-1}})) = \mu(\phi(x_{j_0},\dots,x_{j_{n-1}}))$.
\end{defn}

The previous definition is a natural generalization of the notion of a sequence of indiscernibles. The next proposition (which follows from Lemma 2.10 in \cite{gen_stab_smooth}) is an analogue of the usual characterization of NIP via indiscernibles.

\begin{fact}\label{inv_seq_mes}
Let $T$ be NIP, and $\mu((x_i)_{i<\omega})$ be an indiscernible Keisler measure (over a small model). Let $\phi(x,b)\in L(\U)$. Then the sequence $(\mu(\phi(x_i,b)))_{i<\omega}$ converges as $i$ grows.
\end{fact}

\subsection{The product of types and measures}\label{ssec_prod}

Recall that, given two types $p(x)$ and $q(y)$ with $p$ invariant (over a small set $A$), we can define the product $p\otimes q$ by saying that $p\otimes q\vdash \phi(x,y;c)$ if and only if $p\vdash \phi(x;b,c)$ for some (any) $b\models q_{|Ac}$. It can be checked that this is well defined, and moreover if $q$ is $A$-invariant then $p\otimes q$ is as well. Similarly, if $p$ and $q$ are finitely satisfiable in $M$ (resp. definable over $A$), then $p\otimes q$ too. Additionally, Lemma 2.23 in \cite{NIP_guide} says:

\begin{fact}
Let $p(x)$ be finitely satisfiable in $M$, and $q$ definable over $M$. Then $p_x\otimes q_y=q_y\otimes p_x$.
\end{fact}
 
 Let us now define the product of measures. Let $\mu$ and $\nu$ be two Keisler measures and $T$ NIP. Assume that $\mu$ is $M$-invariant, and let $\phi(x,y;c)\in L(\U)$. Let $N\succ M$ be a model containing $c$, and define the map $f:q\in S_y(N)\mapsto \mu(\phi(x;b,c))\in[0,1]$, where $b$ is any realization of $q$. Then $f$ is Borel (this uses NIP) and we can set \[(\mu\otimes\nu)(\phi(x,y;c))=\int_{S_y(N)}f(q)d\nu(q)\] This defines a new measure $\omega = (\mu\otimes\nu)(x,y)$, satisfying $\omega_{|x}=\mu,\omega_{|y}=\nu$. One can see that, if $p,q$ are types, there product as types coincides with their product as measures.
 
 The product of measures satisfies very similar properties to that of types. For instance, if $\nu$ is $M$-invariant, then $\mu\otimes\nu$ too; if $\mu$ and $\nu$ are finitely satisfiable in $M$ (resp. $M$-definable), then $\mu\otimes\nu$ too (these results are Lemma 1.5 in \cite{gen_stab_smooth}). Note that the product of measures is associative (although this is nontrivial). Also, we have Proposition 7.22 in \cite{NIP_guide}:
 
 \begin{fact}
 Let $\mu(x)$ be finitely satisfiable in $M$ and $\nu(y)$ definable over $M$. Then $\mu_x\otimes\nu_y=\nu_y\otimes\mu_x$.
 \end{fact}

 \begin{remark}
 Let $\mu$ be an $M$-invariant measure and $\nu_1,\nu_2$ having the same restriction to $M$. Then $(\mu\otimes\nu_1)_{|M}=(\mu\otimes\nu_2)_{|M}$.
 
 This follows directly from the definition.
 \end{remark}
 
 Also, the product allows us to form powers. Let $\mu$ be an $M$-invariant measure. We inductively define $\mu^{(n)}(x_0,\dots,x_{n-1})$ by setting $\mu^{(1)}(x_0)=\mu(x_0)$ and $\mu^{(n+1)}(x_0,\dots,x_n)=\mu(x_n)\otimes\mu^{(n)}(x_0,\dots,x_{n-1})$. Passing to the limit, we also define $\mu^{(\omega)}(x_0,x_1,\dots)$.
 
 \begin{fact}\label{Morley_ind}
 The measure $\mu^{(\omega)}$ is indiscernible over $M$.
 \end{fact} 
 
 This last result comes from \cite{gen_stab_smooth} (Lemma 2.14).
 
 \begin{remark}
 Actually, it is possible to define $\mu^{(I)}$ for any linear order $I$. Indeed, if $i_0<i_1<\dots<i_{n-1}$ are in $I$ and $\phi(x_0,\dots,x_{n-1})\in L(\U)$, one can set $\mu^{(I)}(\phi(x_{i_0},\dots,x_{i_{n-1}})) := \mu^{(n)}(\phi(x_0,\dots,x_{n-1}))$. This is clearly well-defined, and indiscernible over $M$ (for instance using indiscernibility of $\mu^{(\omega)}$, which comes itself from a straightforward induction). An interesting property is that, for ordinals, it can still be built by induction: we have $\mu^{(\alpha+1)}(x_0,\dots,x_\alpha) = \mu_{x_\alpha}\otimes\mu^{(\alpha)}$. Indeed, let $\phi(x_0,\dots,x_{n-1})\in L(\U)$, and $i_0<\dots<i_{n-1}$. If $i_{n-1}\neq\alpha$, then 
 
 \[\begin{array}{rcl}
 \mu^{(\alpha+1)}(\phi(x_{i_0},\dots,x_{i_{n-1}})) & = & \mu^{(n)}(\phi(x_0,\dots,x_{n-1})) \\
 & = & \mu^{(\alpha)}(\phi(x_{i_0},\dots,x_{i_{n-1}})) \\
 & = & (\mu_{x_\alpha}\otimes\mu^{(\alpha)})(\phi(x_{i_0},\dots,x_{i_{n-1}}))
 \end{array}\]
 
 On the other hand, if $i_{n-1}=\alpha$, then
 
 \[\begin{array}{rcl}
 \mu^{(\alpha+1)}(\phi(x_{i_0},\dots,x_{i_{n-1}})) & = & \mu^{(n)}(\phi(x_0,\dots,x_{n-1})) \\
 & = & (\mu_{x_{n-1}}\otimes\mu^{(n-1)})(\phi(x_0,\dots,x_{n-1})) \\
 & = & (\mu_{x_\alpha}\otimes\mu^{(\alpha)})(\phi(x_{i_0},\dots,x_{i_{n-1}}))
 \end{array}\]
 
 \end{remark}


\section{Heir property for measure}\label{section_heir}

In the realm of types, the heir property shares a strong connection with definability. More precisely, a type is definable (over $M$) if and only if it has a unique heir (over $M$). This is a useful criterion to show that a type is definable. In particular, heir-coheir duality allows to connect definability with coheirs (this is done by Simon in \cite{dp_types}). In this section, we are going to prove a measure equivalent of it. One of the difficulties is to find a convenient definition for the heir of a measure. Here we will use a definition introduced by Chernikov, Pillay and Simon in \cite{ext_def_NIP_gps} (we state it differently, but equivalence is easy to check). It extends the corresponding notion for types. As in their article, it seems that this definition is more convenient for us than the one proposed in Remark 2.7 of \cite{hrushovski2006groups}.

An analogue of our main result of the section (Proposition \ref{res_heir}) already appears in \cite{hrushovski2006groups} (Remark 2.7). Nonetheless, it was there stated with a weaker definition of an heir, so this new version will better fit our needs.

\begin{defn}
    Let $\mu_x$ be a measure over some set $B$, and $A\subset B$. The measure $\mu$ is said to be an \textit{heir} over $A$ (or to be an heir of $\mu_{|A}$) if, for any $\varepsilon >0, \varphi_i(x;y)\in L(A)$ formulas for $i<n$, and $b\in \U$, there exists some $d\in M$ such that, for every $i$, we have $|\mu(\varphi_i(x;b))-\mu(\varphi_i(x;d))|<\varepsilon$.
\end{defn}

As in the case of types, we can take $d$ to be close to $b$ (i.e. if we choose $\psi(y)\in\tp(b/M)$, $d\in M$ could have been chosen in the previous definition to satisfy $\psi(y)$); this is seen by adding $\psi(y)$ to the $\varphi_i(x;y)$, as soon as $\varepsilon < 1$.

First let us see how to build many heirs. Recall the structure $\tilde M_\mu$ from Subsection 7.1 in \cite{NIP_guide}, which is a way to encode $\mu$ in the language. Let $M$ be an $L$-structure, and $\mu$ a Keisler measure over $M$. Add a new sort for $[0,1]$, endowed with $+$ and $<$. For each $\phi(x;y)\in L$, add a function symbol $f_\phi(y)$, and interpret it by $f_\phi(b) = \mu(\phi(x;b))$. This gives us a new structure $\tilde M_\mu$. If we take an extension $\mathfrak N$ of this structure, we obtain an $L$-structure $N\succ M$, an elementary extension $[0,1]^*\succ [0,1]$ and what we could call a non-standard measure, i.e. something satisfying the axioms of a measure, but with values in $[0,1]^*$. There is a canonical projection $\pi:[0,1]^*\to[0,1]$ (take $\Q$-types and identify with $[0,1]$). Through this projection, we get finally an actual measure $\nu$, extending $\mu$, by setting $\nu(\phi(x;b))=\pi(f_\phi(b))$ for $\phi(x;y)\in L$ and $b\in N$. We say that $\nu$ is the associated (standard) measure to $\mathfrak N$.

\begin{claim}
    Let $\mu$ be a measure over $M$, and $\mathfrak N\succ\tilde M_\mu$ an elementary extension. Then the measure $\nu$ associated to $\mathfrak N$ is an heir of $\mu$.
\end{claim}
\begin{proof}
    Let $\varepsilon>0,\varphi_i(x;y)\in L(M), b\in N$. Denote $r_i:=\nu(\varphi_i(x;b))$. Then \[\tilde N_\nu\models (\exists y)\bigwedge\limits_i |f_{\varphi_i}(y)-r_i|<\varepsilon\]
    Hence $\tilde M_\mu$ also satisfies that formula, proving that $\nu$ is an heir of $\mu$.
\end{proof}

We have seen that, if $\mu$ is an heir over $M$, we could approach the $\mu(\varphi_i(x;b))$ with $d\in M$ arbitrarily close to $b$, in the $M$-definable sense. Conversely, we have:

\begin{lem}\label{exist_heir}
    Let $\mu$ be a measure over $M$. Let $\varphi_i(x;y)\in L(M)$ for $i<n$, and let $b\in \U$. Let also be $r_i\in[0,1]$ for $i<n$. We assume that, for any $\psi(y)\in \tp(b/M)$ and $\varepsilon >0$, there is $d\in M$ such that $d\models\psi(y)$ and $|\mu(\varphi_i(x;d))-r_i|<\varepsilon$ for all $i$. Then there is an heir $\nu$ of $\mu$ such that $\nu(\varphi_i(x;b))=r_i$ for all $i<n$.
\end{lem}
\begin{proof}
    The hypothesis is that the partial type
    \[\pi(x) = \tp(b/M)\cup\{|f_{\varphi_i}(x)-r_i|<\varepsilon:\varepsilon>0,i<n\}\]
    is finitely consistant in $\tilde M_\mu$, hence consistant. Let $\mathfrak N$ be an elementary extension of $\tilde M_\mu$ realizing it. We can expand again $\mathfrak N$ to a monster (whose restriction to $L$ is $\U$ by saturation), and hence obtain $c\in\U$ and $\lambda$ measure over $\U$ such that $c\equiv_M b$ and $\lambda(\varphi_i(x;c))=r_i$ for every $i$. Indeed, $|f_{\varphi_i}(c)-r_i|<\varepsilon$ for all $\varepsilon>0$, so when projecting onto $[0,1]$ we get the desired equality. Moreover, by the above claim, $\lambda$ is an heir of $\mu$. Let $\sigma\in\Aut(\U/M)$ be such that $\sigma(c)=b$. Then $\sigma_*\lambda$ is an heir of $\mu$, and for all $i<n$:\[\sigma_*\lambda(\varphi_i(x;b))=\lambda(\sigma^{-1}(\varphi_i(x;b)))=\lambda(\varphi_i(x;\sigma^{-1}(b))) = \lambda(\varphi_i(x;c))=r_i\]
\end{proof}

We can finally prove the main result of this section:

\begin{prop}\label{res_heir}
    Let $\mu$ be a measure over some small model $M$. It extends to a global $M$-definable measure if and only if it has a unique global heir.
\end{prop}
\begin{proof}
    First assume $\mu$ extends to a global $M$-definable measure $\nu$. Then clearly $\nu$ is an heir, by continuity and density of realized types. Let us now prove that $\nu$ is unique. Let $\lambda$ be a global heir of $\mu$, and let $\varphi(x;y)\in L(M), b\in\U$. For every $\varepsilon>0$, there is $\psi(y)\in\tp(b/M)$ such that $c\models\psi(y)$ implies $|\nu(\varphi(x;c))-\nu(\varphi(x;b))|<\varepsilon$, and then there is $d\in M$ such that $\models \psi(d)$ and $|\lambda(\varphi(x;d))-\lambda(\varphi(x;b))|<\varepsilon$. Then $|\nu(\varphi(x;b))-\lambda((\varphi(x;b))|<2\varepsilon$, and hence $\nu=\lambda$.

    Conversely, let $\nu$ be the only global heir of $\mu$. Note that it is $M$-invariant: if we apply any automorphism that fixes $M$ pointwise, it must also fix $\nu$ by uniqueness of the heir.
    
    \underline{Claim} : Let $\varphi(x;y)\in L(M),b\in\U$ and $\delta>0$. There is $\psi(y)\in\tp(b/M)$ such that, for $d\in M$, if $d\models \psi(y)$, then $|\nu(\varphi(x;d))-\nu(\varphi(x;b))|<\delta$.
    
    Denote $\rho=\nu(\varphi(x;b))$, and let $r\neq \rho$. Apply Lemma \ref{exist_heir} to $\mu$ with $n=1, r_i=r,\varphi_i(x;y)=\varphi(x;y)$ and $b=b$. There is no heir of $\mu$  such that the measure of $\varphi(x;b)$ is $r$ (by uniqueness of the heir of $\mu$). By contraposition we obtain that there exist $\psi_r(y)\in\tp(b/M)$ and $\varepsilon_r>0$ such that, for all $d\in M$ such that $d\models \psi_r(y)$, we have $|\nu(\varphi(x;d))-r|>\varepsilon_r$. Let $\delta >0$. The compact set $[0,1]\setminus ]\rho-\delta,\rho+\delta[$ is covered by the $]r-\varepsilon_r,r+\varepsilon_r[$, hence by a finite number of them, say the $r_i$. Let $\psi(y):=\bigwedge_i\psi_{r_i}(y)\in\tp(b/M)$, and $d\in M$ such that $d\models\psi(y)$. Then $\nu(\varphi(x;d))$ cannot lie in any of the $]r_i-\varepsilon_{r_i},r_i+\varepsilon_{r_i}[$, so $|\nu(\varphi(x;d))-\rho|<\delta$. This proves the claim.

    Let then be $\varphi(x;y),b\in\U,\delta>0$. Let $\psi(y)\in\tp(b/M)$ be given by the claim. Let $c\in\U$ be such that $c\models\psi(y)$. Applying the claim to $c$ and $\varphi(x;y)$, we get $\chi(y)\in\tp(c/M)$ such that, for $d\in \chi(M)$, we have $|\nu(\varphi(x;d))-\nu(\varphi(x;c))|<\delta$. Let $d\in M$ be such that $d\models\psi(y)\wedge\chi(y)$ (it exists, since $c\models \psi(y)\wedge\chi(y)$). Then 
    \[|\nu(\varphi(x;c))-\nu(\varphi(x;b))|\leq|\nu(\varphi(x;c))-\nu(\varphi(x;d))|+|\nu(\varphi(x;d))-\nu(\varphi(x;b))|<2\delta\]
     This shows that $\nu$ is $M$-definable.
\end{proof}



\section{Smooth measures and amalgams}\label{section_smooth}

Smooth measures are measures over a small model which have a unique global extension. They have many useful properties that make them behave as a measure equivalent of realized types. Most of these properties work without NIP, but the key result - the existence, for any measure, of a smooth extension - is a consequence of NIP. Amongst the other properties, a very important one for our purpose is separated amalgamation. When taking two types $p_x$ and $q_y$, we often want to make a single type out of them. If one of them is invariant, this can be done for example via the product, but there is no uniqueness. However, if one of them is a realized type - say $q = \tp(b/\U)$ - there is a unique way of doing so, as there is a unique choice when taking a realization of $q$. Also note that, if we write $p\times q$ this new type, then $(p\times q)_{|A}$ carries the exact same information as $p_{|Ab}$. Separated amalgamation will be a measure equivalent of this construction, and will satisfy similar properties.

\begin{defn}
    A measure $\mu\in\mes_x(M)$ is \textit{smooth} if it has a unique global extension.
\end{defn}

More generally, if $\mu$ is a measure over $A$ (which could be any subset, not necessarily small) and $M\subset A$, then $\mu$ is said to be smooth over $M\subset A$ if $\mu_{|M}$ is smooth. In that case, if $M\prec N\subset A$, $\mu$ is also smooth over $N$. We will sometimes implicitly see smooth measures as global measures (by uniqueness of the extension).

\begin{exemple}
Let $a\in M$. We denote by $\delta_{a,M}$ the measure associated to $\tp(a/M)$, i.e. the measure defined by $\delta_{a,M}(\phi(x)) = 1$ if $a\models\phi(x)$, else $\delta_{a,M}(\phi(x))=0$, for all $\phi(x)\in L(M)$. This measure is smooth: let $\mu$ be a global extension of $\delta_{a,M}$, and $\phi(x)\in L(\U)$. Then if $a\models \phi(x)$ we have $\U\models (a=x)\to\phi(x)$ and $\mu(a=x)=\delta_{a,M}(a=x)=1$, so $\mu(\phi(x))=1$.  If not, then $\U\models\phi(x)\to (a\neq x)$ and $\delta_{a,M}(a\neq x)=0$, so $\mu(\phi(x))=0$. So $\mu$ is unique and coincides with $\delta_{a,\U}$. We will then write $\delta_a$ for $\delta_{a,\U}$. This measure is smooth over any model containing $a$.
\end{exemple}

The following result is Proposition 7.9 in \cite{NIP_guide}.

\begin{fact}
    (T NIP) 
    
    Let $\mu\in\mes_x(M)$. There is $N\succ M$ and a smooth $\mu'\in\mes_x(N)$ extending $\mu$.
\end{fact}

We will often write $\mu'$ for a smooth extension of a measure $\mu$. The following proposition is Lemma 7.8 from \cite{NIP_guide}.

\begin{fact}\label{smooth_psi_theta}
    Let $\mu\in\mes_x(M)$ be smooth, $\phi(x;y)\in L(M)$ and $\varepsilon>0$. Then there are some $\psi_i(y),\theta_i^-(x),\theta_i^+(x)\in L(M), i<n$ such that:
    \begin{enumerate}
        \item The $\psi_i(y)$ form a partition of $M^{|y|}$.
        \item For $i<n$, we have $\models \psi_i(y)\to((\theta_i^-(x)\to\phi(x;y))\wedge (\phi(x;y)\to\theta_i^+(x)))$.
        \item For $i<n$, $\mu(\theta_i^+(x))-\mu(\theta_i^-(x))<\varepsilon$.
    \end{enumerate}
\end{fact}
\begin{remarks}
\begin{itemize}
\item This is an equivalence: if $\mu$ is such that for all $\phi(x;y)\in L(M)$ and $\varepsilon >0$, there are such $\psi_i(y),\theta_i^-(x),\theta_i^+(x)$, then $\mu$ is smooth.
\item If we set $c\in\U,\phi(x;y)\in L(M)$ and $\varepsilon >0$, by choosing the $\psi_i(y)$ in which lies $c$, we obtain two formulas $\theta^-(x),\theta^+(x)\in L(M)$ such that $\theta^-(x)\to \phi(x,c)$, $\phi(x,c)\to \theta^+(x)$ and $\mu(\theta^+(x))-\mu(\theta^+(x))<\varepsilon$. Once more, this is an equivalence: if for any $c\in\U,\phi(x,y)\in L(M)$ and $\varepsilon >0$ we have such $\theta^-(x),\theta^+(x)$, then $\mu$ is smooth. The $\psi_i$ can be found by compactness.
\end{itemize}
\end{remarks}

\begin{cor}\label{petite_base_lisse}
    Let $\mu\in\mes_x(M)$ be smooth. Then it is smooth over a model $N\prec M$ of cardinal at most $|T|$.
\end{cor}
\begin{proof}
Take all the parameters from all the $\psi_i,\theta_i^-$ and $\theta_i^+$ that we get by applying the previous proposition to all formulas $\phi(x;y)\in L$ and all $\varepsilon = \frac 1n$ for $n>0$. This way we get a set $A\subset M$ with $|A|\leq |T|$, which is contained in a model $N\prec M$ with $|N|\leq |T|$.
\end{proof}

\begin{cor}\label{ext_smooth_sat}
    Let $\mu$ be a measure over $M$, and $M\prec^+N$. There is $\mu'$ an extension of $\mu$ which is smooth over $N$.

    Actually, any smooth extension of $\mu$ is an $M$-conjugate of an extension which is smooth over $N$.
\end{cor}
\begin{proof}
    Let $\mu'$ be a smooth (global) extension of $\mu$, and $M_0$ model of cardinality at most $|T|$ such that $\mu'$ is smooth over $M_0$. Let $M_1\subset N$ realize $\tp(M_0/M)$, and $\sigma\in\Aut(\U/M)$ with $\sigma(M_0)=M_1$ (as tuples). Then $\sigma_*\mu'$ is an extension of $\mu$ which is smooth over $N$.
\end{proof}

\begin{defn}\label{def_prod_am}
    Let $\mu_x,\nu_y,\omega_{xy}$ be Keisler measures over $M$. The measure $\omega$ is said to be a \textit{separated amalgam} of $\mu$ and $\nu$ if, for all $\phi(x),\psi(y)\in L(M)$, we have $\omega(\phi(x)\wedge\psi(y))=\mu(\phi(x))\nu(\psi(x))$.
\end{defn}

\begin{remarks}
\begin{itemize}
\item Saying that $\omega$ is a separated amalgam of $\mu$ and $\nu$ means that $\omega$ extends the product measure $\mu\otimes\nu$ on the space $S_x(M)\times S_y(M)$ (here $\otimes$ is the measure-theoretic product, not the model-theoretic one). In particular, such a separated amalgam always exists (by Fact \ref{compact_measures}).
\item A measure $\omega_{xy}$ is said to be an amalgam of $\omega_{|x}$ and $\omega_{|y}$. Sometimes, we just say that $\omega$ is separated for ``$\omega$ is a separated amalgam of $\omega_{|x}$ and $\omega_{|y}$''. This statement is not always true: for instance, if $\mu$ is a measure which is not $0-1$-valued, a nonseparated amalgam of $\mu$ and itself is given by $\omega(\psi(x,y))=\mu(\psi(x,x))$ for all $\psi(x)\in L(M)$.
\item Let $\omega_{xy}$ be a measure, and $\mu :=\omega_{|x}, \nu:=\omega_{|y}$. Assume $\mu$ is a type, i.e. a $0-1$-valued measure. Then the amalgam $\omega$ is separated. Indeed, let $\phi(x),\psi(y)\in L(M)$. If $\mu(\phi(x))=0$, then $\omega(\phi(x)\wedge\psi(y))\leq\omega(\phi(x))=\mu(\phi(x))=0$. Else, we have $\phi(x)=1$. So $\omega(\phi(x)\wedge\psi(y)) = \omega(\psi(y))-\omega((\neg\phi(x))\wedge\psi(x))=\omega(\psi(y))$. In any case, we have $\omega(\phi(x)\wedge\psi(y))=\mu(\phi(x))\nu(\psi(y))$, which is the separation hypothesis.
\end{itemize}
\end{remarks}

The next lemma is Corollary 2.5 in \cite{gen_stab_smooth}, but on any model instead of just in the global setting.

\begin{lem}\label{calcul_prod_smooth}
    Let $\mu_x,\nu_y$ be measures over a model $M$ (which needs not be small). Assume $\mu$ is smooth. Then there is a unique separated amalgam of $\mu$ and $\nu$ (over $M$), which we denote $\mu\times\nu$.
\end{lem}
\begin{proof}
    Let $\phi(x,y)\in L(M)$. Let $\varepsilon>0$, and take $\psi_i(y),\theta_i^-(x),\theta_i^+(x)\in L(M)$ from Proposition \ref{smooth_psi_theta}. Then :
    \[\bigvee_i(\psi_i(y)\wedge\theta_i^-(x))\to\phi(x,y)\text{ and } \phi(x,y)\to\bigvee_i(\psi_i(y)\wedge\theta_i^+(x))\]
    But if $\omega$ is any separated amalgam of $\mu$ and $\nu$, then $\omega(\bigvee_i(\psi_i(y)\wedge\theta_i^-(x)))=\sum_i\mu(\theta_i^-(x))\nu(\psi_i(y))$ and $\omega(\bigvee_i(\psi_i(y,c)\wedge\theta_i^+(x)))=\sum_i\mu(\theta_i^+(x))\nu(\psi_i(y))$, which are $\varepsilon$-close. So there is only one possible value for $\omega(\phi(x,y))$.
    \end{proof}

\begin{remarks}
\begin{itemize}
\item If $\mu,\nu$ are two measures over a model $N$ and $\mu$ is smooth over some $M\prec N$, then we have $(\mu\times\nu)_{|M} = \mu_{|M}\times\nu_{|M}$. In particular, if $\nu_1,\nu_2$ are two measures over $N$ such that $\nu_{1|M}=\nu_{2|M}$, then $(\mu\times\nu_1)_{|M} = (\mu\times\nu_2)_{|M}$.
\item The product $\times$ is commutative. This is straightforward by uniqueness of the separated amalgam.
\item Let $\mu,\nu$ be two global measures. Note that, if $\nu$ is $M$-invariant, then $\nu\otimes\mu$ is a separated amalgam of $\mu$ and $\nu$, so if $\mu$ is smooth $\nu\otimes\mu=\nu\times\mu$. In particular, if $\mu$ is smooth over $M$, then $\mu$ is clearly $M$-invariant; moreover, $\mu_x\otimes\mu_y = \mu_x\times\mu_y=\mu_y\times\mu_x=\mu_y\otimes\mu_x$ so $\mu$ is generically stable (over $M$), i.e. finitely satisfiable in $M$ and definable over $M$ (see Theorem 7.29 in \cite{NIP_guide}).
\item In the case $\mu=\delta_a$ (i.e. $\mu(\phi(x))=1$ if $a\models \phi(x)$, else $0$), the measure $\omega(x,y)$ given by $\omega(\phi(x,y)):=\nu(\phi(a;y))$ is an amalgam of $\delta_a$ and $\nu$. Thus it is separated (since $\delta_a$ is $0-1$-valued), so it is equal to $\delta_a\times\nu$.
\end{itemize}
\end{remarks}

Think of the next result as of the equivalent of the fact that, if $p$ and $q$ are realized types, then $p\times q$ too.

\begin{lem}\label{prod_smooth_smooth}
Let $\mu,\nu$ be smooth measures over $M$. Then $\mu\times\nu$ is smooth over $M$.
\end{lem}
\begin{proof}
Let $\phi(x,y,z)\in L(M)$ and $c\in\U$. Let $\varepsilon>0$, and take $\psi_i(y,z),\theta_i^-(x),\theta_i^+(x)\in L(M),i<n$ from Proposition \ref{smooth_psi_theta}, by smoothness of $\mu$. Then :
    \[\bigvee_i(\psi_i(y,c)\wedge\theta_i^-(x))\to\phi(x,y,c)\text{ and }\phi(x,y,c)\to\bigvee_i(\psi_i(y,c)\wedge\theta_i^+(x))\]
Then, for each $i$, as $\nu$ is smooth over $M$, applying Propositon \ref{smooth_psi_theta} (or more precisely one of the remarks below), we may find some $\chi_i^-(y),\chi_i^+(y)\in L(M)$ such that $\chi_i^-(y)\to\psi_i(y,c)\to\chi_i^+(y)$ and $\nu(\chi_i^+(y))-\nu(\chi_i^-(y))\leq\delta:=\frac \varepsilon n$. We then get:
\[\begin{array}{rcl}
\bigvee_i(\chi_i^-(y)\wedge\theta_i^-(x)) & \to & \bigvee_i(\psi_i(y,c)\wedge\theta_i^-(x)) \\
& \to & \phi(x,y,c) \\
& \text{and } & \\
\phi(x,y,c) & \to & \bigvee_i(\psi_i(y,c)\wedge\theta_i^+(x)) \\
& \to & \bigvee_i(\chi_i^+(y)\wedge\theta_i^+(x))
\end{array}\]

Note that the $\chi_i^-(y)$ are pairwise disjoint, as they are included in the $\psi_i(y,c)$ which partition the $y$-space. This yields in particular $\sum_i\nu(\chi_i^-(y))\leq 1$. Finally:

\[\begin{array}{rcl}
(\mu\times\nu)(\bigvee_i\chi_i^+(y)\wedge\theta_i^+(x)) & \leq & \sum_i (\mu\times\nu)(\chi_i^+(y)\wedge\theta_i^+(x)) \\
& = & \sum_i \nu(\chi_i^+(y))\mu(\theta_i^+(x)) \\
& \leq & \sum_i (\nu(\chi_i^-(y))+\delta)(\mu(\theta_i^-(x))+\varepsilon) \\
& = & \sum_i (\nu\times\mu)(\chi_i^-(y)\wedge\theta_i^-(x)) + \varepsilon \sum_i \nu(\chi_i^-(y)) + \delta \sum_i\mu(\theta_i^-(x)) + n\delta\varepsilon \\ 
& \leq & \sum_i (\nu\times\mu)(\chi_i^-(y)\wedge\theta_i^-(x)) + \varepsilon + \varepsilon + \varepsilon^2 \\
& \leq & (\mu\times\nu)(\bigvee_i\chi_i^-(y)\wedge\theta_i^-(x)) + 3\varepsilon
\end{array}\]

Using the converse of Proposition \ref{smooth_psi_theta} (see the second remark below it), we obtain the smoothness of $\mu\times\nu$.

\end{proof}

\begin{remark}
Let $(\mu_i)_{i\in I}$ be any small family of smooth measures (say $\mu_i$ has variable $x_i$). Then, since a product of two smooth measures is smooth, there is a well-defined product $\prod_{i\in I}\mu_i$, with variable $(x_i)_{i\in I}$. This product is smooth over a small submodel of $\U$. It is the unique separated amalgam of the $\mu_i, i\in I$ (with the natural definition of a separated amalgam in any number of variables).
\end{remark}




The following fact can be found in \cite{NIP_guide} (Lemma 7.17):

\begin{fact}\label{smooth_inv}
    Let $\mu$ be a smooth measure, and assume $\mu$ is $M$-invariant. Then $\mu$ is smooth over $M$.
\end{fact}

We now prove some facts on Boolean algebras of formulas.

\begin{lem}\label{rendre_ex_clos}
Let $\mathcal B$ be a finite Boolean algebra of formulas in variable $xy$. There exists a finite Boolean algebra $\mathcal C$ containing $\mathcal B$ and such that $\mathcal C$ is closed under $\exists y$.
\end{lem}
\begin{proof}
First let $\mathcal B_x$ be the Boolean algebra generated by the formulas $(\exists y)\phi(x,y)$, for $\phi(x,y)\in \mathcal B$. Then define $\mathcal C$ to be the Boolean algebra generated by the $\phi(x,y)\wedge\psi(x)$ for $\phi(x,y)\in \mathcal B$ and $\psi(x)\in \mathcal B_x$. Then $\mathcal C$ is finite and contains $\mathcal B$; it remains to show that it is closed under $\exists y$. Let us find the atoms of $\mathcal C$. First note that \[\neg(\phi(x,y)\wedge\psi(x))  \leftrightarrow  (\neg\phi(x,y)\vee\neg\psi(x))\]
where $\neg\phi(x,y)$ and $\neg\psi(x)$ are both of the form $\phi(x,y)\wedge\psi(x)$. Moreover \[(\phi_1(x,y)\wedge\psi_1(x))\wedge(\phi_2(x,y)\wedge\psi_2(x))  \leftrightarrow  (\phi_1(x,y)\wedge\phi_2(x,y))\wedge(\psi_1(x)\wedge\psi_2(x))\]
which is again of the form $\phi(x,y)\wedge\psi(x)$. So the atoms of $\mathcal C$ are of the form $\phi(x,y)\wedge\psi(x)$. Any formula of $\mathcal C$ can then be written as a finite disjunction $\bigvee_i \phi_i(x,y)\wedge\psi_i(x)$, with $\phi_i(x,y)\in \mathcal B$ and $\psi_i(x)\in\mathcal B_x$. But then 
\[\begin{array}{rcl}
(\exists y)(\bigvee_i \phi_i(x,y)\wedge\psi_i(x)) & \leftrightarrow & \bigvee_i(\exists y)(\phi_i(x,y)\wedge\psi_i(x))\\
& \leftrightarrow &\bigvee_i \psi_i(x)\wedge (\exists y)\phi_i(x,y)
\end{array}\]
Then, $(\exists y)\phi_i(x,y)\in \mathcal B_x$ and $\psi_i(x)\in\mathcal B_x$ by definition, so $\bigvee_i(((\exists y)\phi_i(x,y))\wedge\psi_i(x))\in\mathcal B_x\subset \mathcal C$. So, for $\phi(x,y)\in\mathcal C$, we have $(\exists y)\phi(x,y)\in \mathcal B_x \subset \mathcal C$, i.e. $\mathcal C$ is closed under $\exists y$.
\end{proof}

The following is a useful extension result (similar to ``if $a_1,a_2,b_1\in\U$ and $a_1\equiv_Aa_2$, there is $b_2\in\U$ such that $a_1b_1\equiv_Aa_2b_2$'').

\begin{prop}\label{extension_lem}
    Let $\omega(x,y)$ be a measure over $M$, $\mu'_x$ a smooth extension of $\omega_{|x}$. Then there is $\omega'$ a smooth extension of $\omega$ such that $\omega'_{|x}=\mu'$.
\end{prop}
\begin{proof}
    Let $\mathcal B$ be a finite algebra of $xy$-formulas with parameters in $M$, and $\phi_i(x)\in L(\U), i<n$ a finite number of $x$-formulas. We first want to show that there is a measure $\nu$ over the Boolean algebra generated by $\mathcal B$ and the $\phi_i(x)$ such that $\nu$ coincides with $\omega$ over $\mathcal B$ and $\nu(\phi_i(x))=\mu'(\phi_i(x))$ for all $i$.
    
    Assume without loss of generality that $\mathcal B$ is closed under $\exists y$, using Lemma \ref{rendre_ex_clos}. Denote by $\mathcal B_x$ the Boolean algebra of $x$-formulas of $\mathcal B$. Also assume that the $\phi_i(x)$ are the atoms of a Boolean algebra finer than $\mathcal B_x$. We build $\nu$ by induction on $i$. The induction hypothesis at step $i$ is the following: we have built $\nu_i(x,y)$ measuring the Boolean algebra $\mathcal B_i$ generated by $\mathcal B$ and the $\phi_j(x),j<i$, such that $\nu_i$ coincides with $\omega$ on $\mathcal B$ and with $\mu'$ on the $\phi_j(x),j<i$.
    
    The base step is straightforward. Let $i\geq 0$, and assume we have built $\nu_i$ satisfying the induction hypothesis. Then $\nu_i$ coincides with $\mu'$ over all $x$-formulas. Indeed, let $\mathcal B_{i,x}$ be the Boolean algebra of $x$-formulas in $\mathcal B_i$. It is generated by the $x$-formulas of $\mathcal B$ and the $\phi_j(x),j<i$. Since the $\phi_j(x),j<n$ partition the $x$-space in a finer way than the Boolean algebra of $x$-formulas of $\mathcal B$, the atoms of $\mathcal B_{i,x}$ are of the form $\phi_j(x)$ (for $j<i$) or $\psi(x)\wedge\neg\bigvee_{j\in S}\phi_j(x)$, with $\psi(x)$ an atom of $\mathcal B_x$ and $S = \{j<i: \phi_j(x)\to\psi(x)\}$. Then, $\nu_i(\phi_j(x)) = \mu'(\phi_j(x))$ by induction hypothesis, and
\[\begin{array}{rcl}
\nu_i(\psi(x)\wedge\neg\bigvee_{j\in S}\phi_j(x)) & = & \nu_i(\psi(x))-\sum_{j\in S}\nu_i(\phi_j(x)) \\
& = & \omega(\psi(x))-\sum_{j\in S}\mu'(\phi_j(x)) \\
& = & \mu'(\psi(x))-\sum_{j\in S}\mu'(\phi_j(x)) \\
& = & \mu'(\psi(x)\wedge\neg\bigvee_{j\in S}\phi_j(x))
\end{array}\]
so $\nu_i$ indeed coincides with $\mu'$ on all formulas over which they are both defined.

Note that $\mathcal B_i$ is still closed under $\exists y$. We try and extend $\nu_i$ to $\phi_i(x)$. If some formula $\psi(x,y)\in\mathcal B_i$ implies $\phi_i(x)$, then $(\exists y) \psi(x,y)$ also implies $\phi_i(x)$. But by induction $\nu_i(\psi(x,y))\leq \nu_i((\exists y)\psi(x,y)) = \mu'((\exists y)\psi(x,y))\leq\mu'(\phi_i(x))$. Hence (as we can do the same for $\neg\phi_i(x)$) we can extend $\nu_i$ to a $\nu_{i+1}$ measuring $\phi_i(x)$ with $\nu_{i+1}(\phi_i(x))=\mu'(\phi_i(x))$, using Fact \ref{compact_measures}. Then we are done with the induction.

Let $N$ be a model such that $\mu'$ is smooth over $N$. By compactness of $\mes_{xy}(N)$ (Fact \ref{compact_measures}), if $\mu'$ is smooth over $N$, we obtain $\tilde\omega$ a measure over $N$ such that $\tilde\omega_{|x} = \mu'_{|N}$ and $\tilde\omega_{|M}=\omega$. More precisely, for all $\mathcal B$ finite subalgebra of that generated by $xy$-formulas over $M$ and $x$-formulas over $N$, we have seen that the closed set $F_{\mathcal B}:=\{\tilde\omega\in\mes_{xy}(N):\forall \psi(x,y)\in L_{xy}(M)\cap\mathcal B,\tilde\omega(\psi(x,y)) = \omega(\psi(x,y)) \text{ and }\forall \phi(x)\in L_x(N)\cap \mathcal B, \tilde\omega(\phi(x))=\mu'(\phi(x))\}$ is nonempty. Moreover, this family is clearly stable by finite intersections. Thus the intersection of all $F_{\mathcal B}$ is non empty, which is exactly what we wanted. Any smooth extension $\omega'$ of $\tilde\omega$ has the required properties.
\end{proof}

The reader may have noticed that, in the previous lemma, we don't know if $\omega'$ is separated or not. Actually, another question is whether it would be possible to impose the conjugacy class of $\omega'_{|y}$ among smooth extensions of $\omega_{|y}$. Can we impose one of the two previous properties? Or even better:

\begin{question}
Is it true that, if $\omega\in\mes_{xy}(M),\mu'\in\mes_x(\U),\nu'\in\mes_y(\U)$ are such that $\mu',\nu'$ are smooth, $\mu'_{|M}=\omega_{|x},\nu'_{|M}=\omega_{|y}$, and $\omega$ is separated, there is $\sigma\in\Aut(\U/M)$ such that $\mu'\times(\sigma_*\nu')$ is a smooth extension of $\omega$?
\end{question}

We will now see that ``Morley sequences of invariant measures are indiscernible''.

\begin{prop}\label{Morley_smooth}
Let $M$ be a model, $\mu$ an $M$-invariant measure. Let $(\mu'_i)_{i<\omega}$ be a sequence of smooth measures and $(M_i)_{i<\omega}$ a sequence of small submodels of $\U$ such that each $\mu'_i$ is smooth over $M_{i+1}$, and $M = M_0\prec M_1\prec\dots$. Assume that $\mu'_i$ extends $\mu_{|M_i}$ for all $i<\omega$. Then $\nu(x_0,x_1,\dots):=\mu'_0(x_0)\times\mu'_1(x_1)\times\dots$ is an indiscernible measure over $M$.
\end{prop}
\begin{proof}
We show by induction that $(\mu'_{i-1})\times\dots\times\mu'_0)_{|M} = \mu^{(i)}_{|M}$. Base step is straightforward; let $i\geq 0$ and assume that we know $(\mu'_{i-1})\times\dots\times\mu'_0)_{|M} = \mu^{(i)}_{|M}$.

First note that, for all $i<\omega$, the measure $\mu'_{i-1}(x_{i-1})\times\dots\times\mu'_0(x_0)$ is smooth over $M_i$. By the remarks following Lemma \ref{calcul_prod_smooth}, we get that $(\mu'_i\times\mu'_{i-1}\times\dots\times\mu'_0)_{|M_i} = (\mu\times\mu'_{i-1}\times\dots\times\mu'_0)_{|M_i}$. So $(\mu'_i\times\mu'_{i-1}\times\dots\times\mu'_0)_{|M} = (\mu\times\mu'_{i-1}\times\dots\times\mu'_0)_{|M} = (\mu\otimes(\mu'_{i-1}\times\dots\mu'_0))_{|M}$. By induction, as $\mu$ is $M$-invariant, we then obtain $(\mu'_i\times\mu'_{i-1}\times\dots\times\mu'_0)_{|M} = (\mu\otimes\mu^{(i)})_{|M} = \mu^{(i+1)}(x_0,\dots,x_{i-1},x_i)_{|M}$. So $\nu_{|M} = \mu^{(\omega)}_{|M}$. But this is indiscernible over $M$ (and indiscernibility over $M$ only depends on the restriction to $M$ by definition).
\end{proof}

Let us summarize the results on products of which we will make extensive use later:

\begin{prop}\label{prop_prod}
\begin{itemize}
\item There is a product $\times$ on the space of measures over $M$ (any model) defined a soon as one of the factors is smooth. This product is associative and commutative. If both factors are smooth, it is also smooth (over the same base).
\item There is another product $\otimes$ on the space of global measures defined as soon as one the left factor is invariant. This product is associative. If both factors are invariant, then their product is as well.
\item When $\mu\times\nu$ and $\mu\otimes\nu$ are both defined, they are equal.
\item If $\mu$ is smooth over $M$ and $\nu_{1|M}=\nu_{2|M}$, then $(\mu\times\nu_1)_{|M}=(\mu\times\nu_2)_{|M}$. Actually, if $M\subset A$ and $\nu_{1|A}=\nu_{2|A}$, then $(\mu\times\nu_1)_{|A}=(\mu\times\nu_2)_{|A}$.
\item If $\mu$ is global $M$-invariant and $\nu_{1|M}=\nu_{2|M}$, then $(\mu\otimes\nu_1)_{|M}=(\mu\otimes\nu_2)_{|M}$.
\end{itemize}
\end{prop}
\begin{proof}
Only the `actually' of the fourth point has not been proven yet. It can be seen by looking at how we built the unique separated amalgam. Let $\phi(x,y;c)\in L(A)$ and $\varepsilon > 0$. Let also be $\psi_i(y,z),\theta^-_i(x),\theta_i^+(x)\in L(M)$ associated to the formula $\phi(x;y,z)$ coming from the smoothness of $\mu$ and Proposition \ref{smooth_psi_theta}. Then, we have
\[\begin{array}{rcl}
(\mu\times\nu_1)(\bigvee_i(\psi_i(y,c)\wedge\theta^-_i(x))) & = & \sum_i \mu(\theta_i^-(x))\nu_1(\psi_i(y,c)) \\
& = & \sum_i \mu(\theta_i^-(x))\nu_2(\psi_i(y,c)) \\
& = & (\mu\times\nu_2)(\bigvee_i(\psi_i(y,c)\wedge\theta^-_i(x)))
\end{array}\]
and the same holds for $\theta^+$. Hence we get $(\mu\times\nu_1)(\phi(x,y;c))=(\mu\times\nu_2)(\phi(x,y;c))$.
\end{proof}


\section{An alternative construction for the canonical retraction}\label{section_constr}

In this section, we develop some tools for the study of measures in a locally NIP context. We will more particularly focus on the case where $T$ is NIP, and we add a unary predicate for a submodel. This will yield an interesting application: a construction of the canonical retraction over the space of measures mimicking the type case. This construction, apart from the symmetry, allows us to translate more results from types to measures, and especially one which is crucial for the main theorem.

Moreover, this construction also has the interest of showing how Theorem \ref{loc_smooth_pair} works, and how it makes it possible to use smooth measures jointly with pairs of models - in a context which is not globally NIP anymore.

\subsection{Local smooth measures}
Here the theory $T$ is arbitrary, not necessarily NIP. Let $\Delta$ be a Boolean algebra of formulas. Assume $\Delta$ is NIP (in the sense that all of its formulas are), and write $\Delta(A)$ for the set of formulas of $\Delta$ with $x$ as variable and parameters in $A$ (in this subsection and the next, the variable of our measures will be $x$). Say a $\Delta$-measure is smooth if it has only one global extension (as a $\Delta$-measure).

The following result is the first that leads us to believe that smooth measures could be interesting also in a locally NIP context. Its demonstration is the same as in the global NIP case (see Proposition 7.9 in \cite{NIP_guide}), but we include it for the sake of completeness.

\begin{prop}
Let $\mu$ be a $\Delta(A)$-measure. It has a smooth extension.
\end{prop}
\begin{proof}
It not, for every $\mu_i$ a $\Delta(A_i)$-measure extending $\mu$, there are $\mu_i^0$ and $\mu_i^1$ two distinct extensions of $\mu_i$. Let us build a sequence by setting $\mu_0=\mu,A_0=A$, taking union at limits and by $\mu_{i+1}:=\frac 12(\mu_i^0+\mu_i^1)$ and $A_{i+1}\supset A_i$ such that $\mu^0_{i|A_{i+1}}\neq\mu^1_{i|A_{i+1}}$. Let also be $\phi_i(x;b_i)\in L(A_{i+1})$ and $n_i>0$ such that $|\mu_i^0(\phi_i(x;b_i))-\mu_i^1(\phi_i(x;b_i))|\geq\frac 1{n_i}$. We then obtain such a sequence $(\mu_i,A_i)_{i<|T|^+}$.

For $i<|T|^+$, let $\psi(x)\in L(A_i)$. Then, since $\mu_i^0(\psi(x))) = \mu_i(\psi(x)) = \mu_i^1(\psi(x))$:

\[\begin{array}{rcl}
\mu_{i+1}(\psi(x)\Delta\phi_i(x;b_i)) & = & \frac12(\mu_i^0(\psi(x)\Delta\phi_i(x;b_i))+\mu_i^1(\psi(x)\Delta\phi_i(x;b_i))) \\
& \geq & \frac12(|\mu_i^0(\psi(x))-\mu_i^0(\phi_i(x;b_i))|+|\mu_i^1(\psi(x))-\mu_i^1(\phi_i(x;b_i))|) \\
& \geq & \frac 12|\mu_i^0(\phi_i(x;b_i))-\mu_i^1(\phi_i(x;b_i))| \\
& \geq & \frac 1{2n_i}
\end{array}\]

By the pigeonhole principle, one may assume that $\phi_i(x;y)$ and $n_i$ do not depend on $i$. Let $\nu$ be the union of all $\mu_i,i<|T|^+$. Then, if $i<j<|T|^+$:
\[\nu(\phi(x;b_i)\Delta\phi(x;b_j)) = \mu_{j+1}(\phi(x;b_i)\Delta\phi(x;b_j)) \geq\frac 1{2n}\]
which is a contradiction by NIP (using Lemma 7.6 from \cite{NIP_guide}; this lemma is stated in a globally NIP context, but clearly works as soon as the considered formula is NIP).
\end{proof}

We also have the following result, immediate by compactness but of which we will make extensive use:

\begin{lem}\label{theta_loc}
Let $\mu$ be a global $\Delta$-measure, smooth over $M$. Let $\phi(x)\in \Delta(\U)$ and $\varepsilon >0$. Then there are $\theta^-(x),\theta^+(x)\in \Delta(M)$ such that $\theta^-(x)\to \phi(x)$ and $\phi(x)\to \theta^+(x)$, and $\mu(\theta^+(x))-\mu(\theta^-(x)) < \varepsilon$.
\end{lem}
\begin{proof}
Let $r = \mu(\phi(x))$. Assume for instance that there is no $\theta^-(x)\in \Delta(M)$ such that $\theta^-(x)\to \phi(x)$ and $\mu(\theta^-(x))>r-\frac\varepsilon 2$ (the same can be done for $\theta^+$). Then, by Fact \ref{compact_measures}, there is a global ($L$-)measure $\nu$ extending $\mu_{|M}$ such that $\nu(\phi(x))\leq r -\frac\varepsilon 2$. But we have $\nu_{|\Delta(\U)}\neq\mu$, which contradicts the smoothness of $\mu_{|M}$.
\end{proof}

Let us recall that a $\Delta$-measure $\mu$ is generically stable over $M$ if it is definable and finitely satisfiable in $M$ (with the same definitions of these notions as the ones given in the introduction, except that we restrict ourselves to $\Delta$-formulas).

\begin{prop}\label{loc_gen_stable}
If a $\Delta$-measure $\mu$ is smooth over $M$, it is generically stable over $M$.
\end{prop}
\begin{proof}
We need to check it is finitely satisfiable in $M$ and definable over $M$.

Let $\phi(x;b)\in \Delta(\U)$ be such that $\mu(\phi(x;b))>0$. Then, by Lemma \ref{theta_loc}, there is $\theta^-(x)\in\Delta(M)$ such that $\theta^-(x)\to\phi(x;b)$ and $\mu(\theta^-(x))>0$. But then $\theta^-(x)$ is satisfiable in $M$, so $\phi(x;b)$ too.

For definability, let $r>0$ and assume that $\mu(\phi(x;b))<r$ for some $b\in\U$. Let $\varepsilon>0$ be such that $\mu(\phi(x;b))<r-\varepsilon$. By Lemma \ref{theta_loc}, there is $\theta^+(x)\in L(M)$ such that $\mu(\theta^+(x))<r$ and $\phi(x;b)\to\theta^+(x)$. Then, there is $\psi(y)\in \Delta(M)$ such that $c\models\psi(y)$ if and only if $\phi(x;c)\to\theta^+(x)$ holds. In particular, $b\models\psi(y)$. So for any $c\in\U$ such that $c\models\psi(y)$ (i.e. for any $c$ in an neighbourhood of $b$) we have $\mu(\phi(x;c))<r$. So we have shown that $\{b\in\U:\mu(\phi(x;b))<r\}$ is open for the $M$-logic topology, i.e. $\mu$ is definable.
\end{proof}

\subsection{Smooth measures and pairs of models}
Let $M$ be a model of an NIP theory, and $A\subset M$. We take a unary predicate $P$ for $A$ in $M$, and denote by $L_P$ the language obtained by adding $P$ to $L$. We denote by $\Delta$ the Boolean algebra generated by $L$-formulas and $P(x)$, and by $(M,A)$ the $L_P$-structure obtained by setting $P(M)=A$. Also, if $x=x_1\dots x_n$ is an $n$-ary variable, we set $P(x) = P(x_1)\wedge\dots\wedge P(x_n)$. Let us start with the following observations:

\begin{claim}
	\begin{itemize}
		\item If $\mu$ is an $L$-measure on $(M,A)$ (or on a subalgebra), it can be extended to a $\Delta$-measure by setting $\mu(P(x))=1$ if and only if $\mu$ is finitely satisfiable in $A$.
		\item If $T$ is NIP, then all formulas in $\Delta$ are NIP. Notably, if $\mu$ is a $\Delta$-measure on a small set $B$, it can be extended to a smooth $\Delta$-measure.
		\item The formulas of $\Delta(M)$ are of exactly those of the form $(\psi(x)\wedge P(x)) \vee (\chi(x)\wedge\neg P(x))$, with $\psi(x),\chi(x)\in L(M)$.
	\end{itemize}
\end{claim}

The first observation directly comes from Fact \ref{compact_measures}. The second one is an immediate consequence of the fact that $P(x)$ is NIP (since it has no parameter) and NIP formulas are preserved under Boolean combinations. The third one is straightforward.


Now assume that $A$ is an elementary submodel of $M$. We change our notations from $(M,A)$ to $(M',M)$. Our goal will be to prove the measure analogue of the following trivial fact: let $p$ be a type over $M'$, finitely satisfiable in $M$. Assume that $p\cup\{P(x)\}$ is a realized type. Then $p$ is realized in $M$. Although this is completely straightforward when speaking of types, the measure analogue is quite difficult - and requires NIP.

\begin{thm}\label{loc_smooth_pair}
Let $\mu$ be an $L$-measure over $M'$, finitely satisfiable in $M$. Assume that the $\Delta$-measure $\tilde\mu$ defined by setting $\tilde\mu(P(x))=1$ and $\tilde\mu_{|L}=\mu$ is smooth (over $M'$). Then $\mu$ is smooth over $M$.
\end{thm}
\begin{proof}
Expand $\U$ to $(\U,P(\U))\succ (M',M)$ - this is possible by saturation. 
We still denote by $\tilde\mu$ the unique global extension of $\tilde\mu$, and we also denote by $\mu$ the global measure defined by $\mu:=\tilde\mu_{|L}$. This is coherent with our previous notation on $M$.

\underline{Claim:} Let $(M',M)\prec^+(N',N)\prec (N'_1,N_1)$, and $\nu$ be such that $\nu_{|N} = \mu_{|N}$. Then $\nu_{|N_1}=\mu_{|N_1}$.

Let $\phi(x)\in L(N_1)$ and $\varepsilon >0$. By Lemma \ref{theta_loc}, there are $\Delta$-formulas $\phi^-(x),\phi^+(x)\in \Delta(M')$ such that
\[\phi^-(x)\to\phi(x)\text{ and }\phi(x)\to\phi^+(x)\]
and also $\tilde\mu(\phi^+(x))-\tilde\mu(\phi^-(x))\leq\varepsilon$. We can write $\phi^-(x) = (\theta^-(x)\wedge P(x))\vee(\chi^-(x)\wedge\neg P(x))$, with $\theta^-(x),\chi^-(x)\in L(M')$. But then $\theta^-(x)\wedge P(x)\to\phi^-(x)$, and $\tilde\mu(\theta^-(x)\wedge P(x)) = \tilde\mu(\phi^-(x))$ because $\tilde\mu(P(x))=1$. So we can replace $\phi^-(x)$ by $\theta^-(x)\wedge P(x)$. Similarly, we can replace $\phi^+(x)$ by $\theta^+(x)\vee\neg P(x)$ for some $\theta^+(x)\in L(M')$. So we have got:
\[\theta^-(x)\wedge P(x)\to \phi(x)\text{ and }\phi(x)\to \theta^+(x)\vee\neg P(x)\]
with $\theta^-(x),\theta^+(x)\in L(M')$ and $\tilde\mu(\theta^+(x))-\tilde\mu(\theta^-(x)) \leq\varepsilon$ (or equivalently $\mu(\theta^+(x))-\mu(\theta^-(x))\leq\varepsilon$). Taking conjuction with $P(x)$ changes the equation to:

\begin{equation}\label{theta_phi}
\theta^-(x)\wedge P(x) \to \phi(x)\wedge P(x) \text{ and }\phi(x)\wedge P(x)\to \theta^+(x)\wedge P(x)
\end{equation}

Now, by Theorem 3.13 in \cite{NIP_guide}, since $M\prec^+ N$, NIP ensures that there are $\psi^-(x),\psi^+(x)\in L(N)$ such that $\psi^-(M)=\theta^-(M),\psi^+(M)=\theta^+(M)$ and $\psi^-(N)\subset \theta^-(N),\psi^+(N)\supset\theta^+(N)$. Then, the inclusions $\psi^-(N)\subset\theta^-(N)$ and $\theta^+(N)\subset\psi^+(N)$ can be changed to $(\psi^-\wedge P)(N')\subset (\theta^-\wedge P)(N')$ and $(\theta^+\wedge P)(N')\subset (\psi^+\wedge P)(N')$. But all of $\psi^-\wedge P,\psi^+\wedge P,\theta^-\wedge P,\theta^+\wedge P$ are $L_P$-formulas with parameters in $N'$. So these inclusions mean:
\[N'\models \psi^-(x)\wedge P(x)\to\theta^-(x)\wedge P(x)\text{ and } N'\models \theta^+(x)\wedge P(x)\to\psi^+(x)\wedge P(x)\]

But $N'\prec N'_1$, so this still holds in $N'_1$. Then, combining with Equation \ref{theta_phi}:

\[N'_1\models \psi^-(x)\wedge P(x)\to\phi(x)\wedge P(x)\text{ and }N'_1\models\phi(x)\wedge P(x)\to \psi^+(x)\wedge P(x)\]

Hence $\psi^-(N_1)\subset\phi(N_1)\subset\psi^+(N_1)$. All of these are $L$-formulas with parameters in $N_1$. So
\[N_1\models \psi^-(x)\to \phi(x)\text{ and } N_1\models \phi(x)\to\psi^+(x)\]
But $N_1\prec N'_1$, so this still holds in $N'_1$. Then, we get
\[\begin{array}{rcccl}
\mu(\psi^-(x))&\leq&\mu(\phi(x))&\leq&\mu(\psi^+(x))\\
\nu(\psi^-(x))&\leq&\nu(\phi(x))&\leq&\nu(\psi^+(x))
\end{array}\] But now, recall that $\psi^-(M)=\theta^-(M),\psi^+(M)=\theta^+(M)$ and $\mu$ is finitely satisfiable in $M$. So $\mu(\psi^-(x)) = \mu(\theta^-(x)), \mu(\psi^+(x))=\mu(\theta^+(x))$, which are $\varepsilon$-close. In addition, $\mu$ and $\nu$ coincide over $N$ by assumption, so $\nu(\psi^-(x))=\mu(\psi^-(x))$ and $\nu(\psi^+(x))=\mu(\psi^+(x))$. Finally, $\mu(\phi(x))$ and $\nu(\phi(x))$ are in the same interval $[\mu(\theta^-(x)),\mu(\theta^+(x))]$ of length at most $\varepsilon$, and this holds for all $\varepsilon >0$. So $\mu(\phi(x))=\nu(\phi(x))$, concluding the proof of this claim.

Set $(M',M)\prec^+(N',N)\prec^+ (N_2',N_2)$. Let $\nu$ be a smooth extension of $\mu_{|N}$. Let $N_1\prec\U$ be such that $\nu$ is smooth over $N_1$ and $|N_1|\leq |T|$ (it exists by Corollary \ref{petite_base_lisse}). We fix an enumeration $(n_i)_{i\in I}$ of $N_1$, and write $y_i$ for the variable associated to $n_i$. The partial $L_P$-type $\tp_L((n_i)_{\in I}/N)\cup\bigcup_{i\in I}\{P(y_i)\}$ is consistent (since any formula of $\tp_L((n_i)_{i\in I}/N)$ is realized in $N$, which is a model in the language $L$). Let $(n'_i)_{i\in I}\subset (N'_2,N_2)$ be a realization of it. We may write $\hat N_1:=\{n'_i|i\in I\}$, which is in natural bijection with $N_1$. We then have $N\prec \hat N_1\prec N_2$. Let $\sigma\in \Aut_L(\U/N)$ be such that $\sigma(n_i)=n'_i$ for all $i\in I$, and $\hat\nu :=\sigma_*\nu$. We have $\hat\nu_{|N} = \nu_{|N}=\mu_{|N}$, so the claim shows that $\hat\nu_{|N_2}= \mu_{|N_2}$. But $\hat\nu$ is smooth over $N_2$, so actually $\mu$ is smooth over $N_2$.

It only remains to see that $\mu$ is $M$-invariant. Let $\phi(x)\in L(\U)$ be such that $\mu(\phi(x))>0$. Then we also have $\tilde\mu(\phi(x))>0$, and since $\tilde\mu(P(x))=1$ we get $\tilde\mu(\phi(x)\wedge P(x))>0$. But $\tilde\mu$ is smooth over $M'$, hence (by Lemma \ref{loc_gen_stable}) finitely satisfiable in $M'$. So the formula $\phi(x)\wedge P(x)$ is satisified in $M'$, which means that $\phi(x)$ is satisified in $M$.

We have seen that $\mu$ is finitely satisfiable in $M$, and in particular it is $M$-invariant. We find that $\mu$ is smooth and $M$-invariant, so it is actually smooth over $M$ by Fact \ref{smooth_inv}.

\end{proof}

This result seems quite important: it allows to use smooth measures on pairs of models. As pairs of models are a powerful tool in model theory, this theorem certainly makes it possible to apply the powerful language of smooth measures in many new contexts. The next subsection will be an example of it.
  
\subsection{An alternative construction for the canonical retraction}\label{ssec_alt_cons}
Let $T$ be NIP. The canonical retraction for measures is defined in \cite{ext_def_NIP_gps} by pushing forward the one for types. In this section, we propose a construction which mimics the one for types, and show that it coincides with the usual one. The big advantage of this is that it allows to translate more results directly from types to measures, as we will see with Corollary \ref{F_M_comm_prod}.

Fix $\mu$ an $M$-invariant measure. Let $M\prec^+M'$ and $(M',M)\prec^+(N',N)$ be elementary extensions. Let us show that, if $\mu_1,\mu_2$ are two global $\Delta$-measures extending $\mu_{|N}$ with $\mu_1(P(x))=\mu_2(P(x))=1$, then $\mu_1$ and $\mu_2$ have the same $L$-reduct to $M'$. If not, take two such $\mu_1,\mu_2$ and $\phi(x;b)\in L(M')$ such that $\mu_1(\phi(x;b))\neq\mu_2(\phi(x;b))$. Let us build a sequence as follows: set $M_0=M$, and then, given a submodel $M\prec M_i\prec N$, with $|M_i|\leq |M|+|T|$, if $i$ is even, we take $\nu_i$ to be an extension of $\mu_{1|M_ib}$ which is smooth over $N'$, and if $i$ is odd, of $\mu_{2|M_ib}$. We have $\nu_i(P(x))=1$, so the $L$-reduct of $\nu_i$ is smooth over $N$. Then we take $M_{i+1}\prec N$ containing $M_i$ (elementarily) and such that $\nu_i$ is smooth over $M_{i+1}$, with $|M_{i+1}|\leq |M|+|T|$ (this is possible by Corollary \ref{petite_base_lisse}). Then the $L$-measure $\prod_{i<\omega}\nu_{i|L}$ (recall from the remark below Lemma \ref{prod_smooth_smooth} that this is well defined) is indiscernible over $M$ (by Proposition \ref{Morley_smooth}). But $\nu_{2i}(\phi(x;b))=\mu_1(\phi(x;b))$ and $\nu_{2i+1}(\phi(x;b))=\mu_2(\phi(x;b))$, which contradicts NIP (or rather Fact \ref{inv_seq_mes}).

We can then set $\tilde F_M(\mu)$ to be the only $M$-invariant measure such that $\tilde F_M(\mu)_{|L,M'} = \nu_{|L,M'}$, where $\nu$ is a $\Delta$-measure extending $\mu_{|N}$ with $\nu(P(x))=1$. This is clearly well-defined and finitely satisfiable in $M$ (which can be proven exactly in the same way as for the case of types). We also have $\tilde F_M(\mu)_{|M}=\mu_{|M}$. In the case where $\mu=p$ is an $M$-invariant type, the construction makes it clear that $\tilde F_M(p)=F_M(p)$. The map $\tilde F_M$ is clearly affine. If $\mu$ is finitely satisfiable in $M$, then it is possible to extend $\mu_{|N'}$ by setting $\mu(P(x))=1$, so by definition $\tilde F_M(\mu)=\mu$.

\begin{lem}
$\tilde F_M$ is continuous.
\end{lem}
\begin{proof}
Let $\phi(x;b)\in L(\U)$, and $r\in [0,1]$. Take $M'$ to contain $b$. If $\tilde F_M(\mu)(\phi(x;b))> r$, then $\nu(\phi(x;b))>r$ for all $\nu$ being a $\Delta$-measure extending $\mu_{|N}$ with $\nu(P(x))=1$. Then, by Fact \ref{compact_measures}, there is $\psi(x;c)\in L(N)$ such that $\mu(\psi(x;c))> r$ and $\psi(x;c)\wedge P(x)\to\phi(x;b)$. Hence,  if $\lambda$ is another $M$-invariant measure with $\lambda(\psi(x;c))>r$, we have $\tilde F_M(\lambda)(\phi(x;b)) >r$. In other words, $\tilde F_M(\lambda)(\phi(x;b))>r$ for all $M$-invariant measures in an open neighbourhood of $\mu$, i.e. $\tilde F_M$ is continuous.
\end{proof}

The following lemma is directly adapted from Proposition 7.11 in \cite{NIP_guide}.

\begin{lem}\label{approx_mes}
Let $\mu\in\mes_x(\U)$ be a global Keisler measure. Let $\phi_i(x;y)\in L,i<n$ be formulas, and let $\varepsilon >0$. Then there exist types $p_0,\dots,p_{N-1}$ in the support of $\mu$ such that, for all $b\in\U$ and $i<n$, we have:
\[|\mu(\phi(x;b))-\frac1N\sum_{k<N}p_k(\phi(x;b))|<\varepsilon\]
\end{lem}
\begin{proof}
This is an immediate consequence of Proposition 7.11 in \cite{NIP_guide} and Lemma 3.12 in Chapter 1 (An introduction to stability theory) of \cite{munster}. Indeed, Proposition 7.11 says that this result is true for one formula, and then Lemma 3.12 allows to go directly from the one-formula version to the $n$-formulas version.
\end{proof}

\begin{thm}\label{F_M_F_M}
$\tilde F_M = F_M$
\end{thm}
\begin{proof}
We know they coincide over $M$-invariant types. They are both affine, hence coincide over the convex hull of $M$-invariant types. They are both continuous, hence coincide over the closure of the convex hull of $M$-invariant types. It remains to show that this is the whole space of $M$-invariant measures. But this follows from the previous lemma. Indeed, let $\mu$ be an $M$-invariant measure, and let $\mu\in V$ be an open neighbourhood of $\mu$. There are $\phi_i(x;y),i<n$ formulas, $b\in\U$ and $\varepsilon >0$ such that any $M$-invariant measure $\nu$ with $|\nu(\phi_i(x;b))-\mu(\phi_i(x;b))|<\varepsilon$ for all $i$ lies in $V$. Then, by Lemma \ref{approx_mes}, there are types $p_0,\dots,p_{N-1}$ in the support of $\mu$ such that $|\frac 1N\sum_{k<N}p_k(\phi_i(x;b))-\mu(\phi_i(x;b))|<\varepsilon$ for all $i$. But then the measure $\frac 1N\sum_{k<N}p_k$ lies in $V$. Moreover, the types $p_k$ are $M$-invariant, since they are in the support of an $M$-invariant measure and $T$ is NIP (see below Definition 7.16 in \cite{NIP_guide}). So $V$ has a nonempty intersection with the convex hull of $M$-invariant types. Hence $\mu$ belongs to the closed convex hull of $M$-invariant types.
\end{proof}

Finally, we give an application of this new construction:

\begin{cor}\label{F_M_comm_prod}
Let $\mu,\nu$ be $M$-invariant measures, and assume $\nu$ is finitely satisfiable in $M$. Then $F_M(\nu\otimes\mu)=\nu\otimes F_M(\mu)$.
\end{cor}
\begin{proof}
Let $M\prec^+M'$, and let $(M',M)\prec^+(N',N)\prec^+(N_1',N_1)\prec^+(N'_2,N_2)$. Then $N\prec^+N_1$ so, by Corollary \ref{ext_smooth_sat}, there is $\mu'$ smooth over $N_1$ extending $\mu_{|N}$. The measure $\nu$ is finitely satisfiable in $M$, so we can extend $\nu$ to a global $\Delta$-measure $\lambda$ by setting $\lambda(P(x))=1$ (and $\lambda_{|L}=\nu$). Then take $\lambda'$ an extension of the $\Delta$-measure $\lambda_{|N'_1}$ which is smooth over $N'_2$ (again, it exists by Corollary \ref{ext_smooth_sat}). The $L$-measure $\nu':= \lambda'_{|L}$ is smooth over $N_2$ by Theorem \ref{loc_smooth_pair}, as $\lambda'(P(x))=1$. Additionnally, $\nu'_{|N'_1}=\nu_{|N'_1}$ by construction. So $(\nu'\times\mu')_{|N_1} = (\nu\times\mu')_{|N_1}$ since $\nu'_{|N_1}=\nu_{|N_1}$ and $\mu'$ is smooth over $N_1$. Then, \[(\nu'\times\mu')_{|N} = (\nu\times\mu')_{|N} = (\nu\otimes\mu)_{|N}\] where the last equality comes from $\mu'_{|N} = \mu_{|N}$ and $N$-invariance of $\nu$.

Moreover, $\nu'\times\mu'$ is smooth over $N_2$ by Lemma \ref{prod_smooth_smooth}, so in particular it is finitely satisfiable in $N_2$ (by Lemma \ref{loc_gen_stable}). When working in the pair $(N'_2,N_2)$, this means that $\nu'\times\mu'$ can be extended by setting that the measure of $P(x,y) = P(x)\wedge P(y)$ is $1$. Hence, by construction of $\tilde F_M$ and Theorem \ref{F_M_F_M}, as $(M',M)\prec^+(N'_2,N_2)$, we have: \[F_M(\nu\otimes\mu)_{|M'} = (\nu'\times\mu')_{|M'}\]

On the other hand, since $\mu'$ is smooth over $N'_1$ and $\nu'_{|N'_1} = \nu_{|N'_1}$, we have $(\nu'\times\mu')_{|N'_1} = (\nu\times\mu')_{|N'_1}$. In particular, $(\nu'\times\mu')_{|M'} = (\nu\times\mu)_{|M'}$. But $\nu$ is invariant over $M'$ (it is finitely satisfiable in $M$) and $\mu'_{|M'} = F_M(\mu)_{|M'}$ by Theorem \ref{F_M_F_M}. So $(\nu\otimes\mu')_{|M'} = (\nu\otimes F_M(\mu))_{|M'}$, and we get $(\nu'\times\mu')_{|M'} = (\nu\times\mu')_{|M'} = (\nu\otimes\mu')_{|M'} = (\nu\otimes F_M(\mu))_{|M'}$. Combining our two equalities yields $F_M(\nu\otimes\mu)_{|M'} =(\nu'\times\mu')_{|M'} = (\nu\otimes F_M(\mu))_{|M'}$. By saturation, we get the desired equality.
\end{proof}


\section{Characterization of definable measures}\label{section_def}

In this section, we will adapt to measures a few results already known for types, in particular Fact \ref{main_thm_type}. Smooth measures will prove to be a very powerful tool for this purpose. The first subsection gives an equivalent of the heir-coheir duality and applies it to show that being a ($M$-)definable measure is the same as ($M$-)commuting to all finitely satisifiable measures (in $M$). 
In the third subsection, we prove our theorem. In that subsection, we follow the proof of \cite{inv_type_NIP}, making an extensive use of our preceding results on smooth measures, illustrating their handyness.

The theory $T$ is NIP in all the section.

\subsection{Definability and commutation}

We begin by an equivalent of the heir-coheir duality:

\begin{prop}
    Let $\mu$ be a global heir over $M$, $b\in\U$, and $\mu'$ a smooth extension of $\mu$ over $N$ a model containing $Mb$. There exists a type $q$ coheir of $\tp(b/M)$ such that $(\delta_b\times \mu')_{|M}=(q\times\mu')_{|M}$.
\end{prop}
\begin{proof}
    \underline{Claim:} there exists a global measure $\omega$ such that:
    \begin{enumerate}
        \item $\omega_{|M}=(\mu'\times \delta_b)_{|M}$
        \item $\omega_{|x}=\mu'$
        \item $\omega_{|y}$ is a type which is a coheir over $M$.
    \end{enumerate}
    By compactness of $\mes_{xy}(\U)$, it is enough to show the following finite version:
    
    \underline{Claim 2:} Let $\phi_i(x)\in L(\U),\psi_j(x,y)\in L(M), \chi_l(y)\in L(\U)$. Assume the $\phi_i$, the $\psi_j$ and the $\chi_l$ all partition their respective space. There exists a measure $\tilde\omega$ over $\mathcal B$ the Boolean algebra generated by the $\phi_i,\psi_j$ and $\chi_l$ such that:
    \begin{enumerate}
    \item For all $j$, we have $\tilde\omega(\psi_j(x,y))=(\mu'\times\delta_b)(\psi_j(x,y))$.
    \item For all $i$, we have $\tilde\omega(\phi_i(x))=\mu'(\phi_i(x))$.
    \item For all $l$, $\tilde\omega(\chi_l(y))$ is equal to $0$ or $1$. Moreover, if it is equal to $1$, then the formula $\chi_l(y)$ is satisfiable in $M$.
    \end{enumerate}
    
    We start by showing the second claim. Let $\varepsilon >0$. By the heir property, there is $d\in M$ such that, for all $j$, we have $|\mu(\psi_j(x;b))-\mu(\psi_j(x;d))|<\varepsilon$. Take $\omega_{\varepsilon}$ to be defined by
    \[\omega_{\varepsilon}(\phi_i(x)\wedge\psi_j(x,y)\wedge\chi_l(y)):=\mu'(\phi_i(x)\wedge\psi_j(x,d)\wedge\chi_l(d))\]
    This trivially defines a measure over $\mathcal B$. Write $a\approx^{\varepsilon}b$ for $|a-b|<\varepsilon$. Then, $\omega_{\varepsilon}$ satisfies:
    \begin{enumerate}
        \item $\omega_{\varepsilon}(\psi_j(x,y))=\mu(\psi_j(x;d))\approx^{\varepsilon}\mu(\psi_j(x;b))=\mu'(\psi_j(x;b))=(\mu'\times \delta_b)(\psi_j(x,y))$
        \item $\omega_{\varepsilon}(\phi_i(x))=\mu'(\phi_i(x))$
        \item $\omega_{\varepsilon}(\chi_l(y))=0\text{ or }1$, and if $\chi_l(M)=\emptyset$, then $\omega_{\varepsilon}(\chi_l(y))=0$.
    \end{enumerate}
    We then take $\varepsilon \to 0$. By compactness of the space of measures over $\mathcal B$ (which is a metrizable space), we can find a convergent subsequence of $(\omega_{\frac 1n})_{n<\omega}$, and denote by $\tilde\omega$ its limit. It satisfies the three hypotheses of the second fact.
    
    Then, compactness of $\mes_{xy}(\U)$ ensures we can find the desired measure $\omega$ (as in the proof of Proposition \ref{extension_lem}).
    
    Let $q$ be the type associated to the measure $\omega_{|y}$ (i.e. $q\vdash \chi(y)\iff \omega(\chi(y))=1$). It is a type, coheir over $M$ by construction. Moreover, $\omega = q\times \mu'$ (by the remarks following Definition \ref{def_prod_am}), concluding the proof.
\end{proof}

Recall that, if $\mu$ is a definable measure and $\nu$ a finitely satisfiable measure, then they commute (see Proposition 7.22 in \cite{NIP_guide}). We will here see a converse to this, similar to the one obtained for types in Lemma 6 of \cite{dp_types}.

\begin{lem}\label{lem_M_comm}
    Let $\mu$ be an $M$-invariant measure. Then it is definable if and only if, for every type $q$ finitely satisfiable in $M$, we have $(\mu\otimes q)_{|M}=(q\otimes\mu)_{|M}$.
\end{lem}
\begin{proof}
    The direct implication is already known, so let us do the converse. By Proposition \ref{res_heir}, it is enough to see that $\mu$ is the unique global heir of $\mu_{|M}$. Let $\nu$ be such an heir, and $\phi(x;b)\in L(\U)$. Let $N$ be a small model containing $Mb$, and $\nu'$ a smooth extension of $\nu_{|N}$. Let $q$ be finitely satisfiable over $M$ be given by the last proposition (applied to $\nu'$). Then:
    \[\begin{array}{rcll}
        \nu(\phi(x;b)) & = & \nu'(\phi(x;b))\\
        & = & (\nu'\times\delta_b)(\phi(x,y))\\
        & = & (\nu'\times q)(\phi(x,y))\\
        & = & (q\times\nu')(\phi(x,y))\\
        & = & (q\otimes\nu')(\phi(x,y))\\
        & = & (q\otimes\nu)(\phi(x,y))\\
        & = & (q\otimes\mu)(\phi(x,y)) & \text{because }q\text{ is }M\text{-invariant and }\mu_{|M}=\nu_{|M}\\
        & = & (\mu\otimes q)(\phi(x,y)) & \text{by assumption}\\
        & = & \mu(\phi(x;b)) & \text{by definition of the product}
    \end{array}\]
    So $\mu = \nu$, which means that $\mu$ is the only global heir of $\mu_{|M}$.
\end{proof}

\subsection{The canonical retraction}
This subsection only contains the proof of our main result, which is quite long. It is adapted from Proposition 3.8 in \cite{inv_types_NIP}, but with some changes; in particular, the use of reverse types (which is fundamental is Simon's proof) has been completely avoided here - although it is possible to define a notion of reverse measure analogous to that of reverse type, and which has similar properties. The proof involves a lot of manipulation of $\times$ and $\otimes$.

Recall the notion of infinite powers of measures from the definition following Fact \ref{Morley_ind}: if $I$ is any infinite linear order and $\mu_x$ an invariant measure, the measure $\mu^{(I)}$ is defined as the measure with variables $(x_i)_{i\in I}$ such that, for any $i_0<\dots i_{n-1}$ in $I$, we have $\mu^{(I)}_{|x_{i_0}\dots x_{i_{n-1}}} = \mu^{(n)}(x_{i_0},\dots,x_{i_{n-1}})$ (where restriction to variables means that we consider only the measure over these variables; we keep this notation is what will follow). We also introduce the notation $\mu^{(I)}_{rev}$, defined as the measure $\mu^{(I^{opp})}$ where $I^{opp}$ is the opposite order to $I$. Note that, in the same way that we proved that $\mu^{(\alpha+1)} = \mu_{x_\alpha}\otimes\mu^{(\alpha)}$ (just after Fact \ref{Morley_ind}), it holds that $\mu^{(\alpha+1)}_{rev} = \mu^{(\alpha)}_{rev}\otimes \mu_{x_\alpha}$.

The following trivial lemma will be very useful in the proof of the theorem.

\begin{lem}
Let $\mu(x,y)$ be a global $M$-invariant measure. Then $F_M(\mu)_{|x} = F_M(\mu_{|x})$.
\end{lem}
\begin{proof}
Let $M\prec M'^+$ a $(|M|+|T|)^+$-saturated extension, and $(M',M)\prec^+ (N',N)$. Let $\nu$ be a $\Delta$-measure such that $\nu_{|L,N} = \mu_{|N}$ and $\nu(P(x)\wedge P(y))=1$ (where $P$ is the unary predicate for $M$ in $M'$ and $\Delta$ is generated by $L$ and $P$ as in Subsection \ref{ssec_alt_cons}). Then $F_M(\mu)_{|M'} = \nu_{|L,M'}$. Moreover, since $\nu_{|x,L,N} = \mu_{|x,N}$ and $\nu(P(x))=1$, we also have that $F_M(\mu_{|x})_{|M'} = \nu_{|x,L,M'}$. Then $F_M(\mu)_{|x,M'} = F_M(\mu_{|x})_{|M'}$. But $F_M(\mu)_{|x}$ and $F_M(\mu_{|x})$ are both $M$-invariant, so we conclude by saturation.
\end{proof}

We now come to the proof of the main theorem.

\begin{thm}\label{thm_F_M_def}
    Let $\mu$ be an $M$-invariant measure. It is definable if and only if it commutes with $F_M(\mu)$.
\end{thm}
\begin{proof}
    It is known that, if $\mu$ is definable, it commutes with $F_M(\mu)$ (which is finitely satisfiable; this follows from Proposition 7.22 in \cite{NIP_guide}).

    Conversely, assume that $\mu$ commutes with $F_M(\mu)$ but is not definable. Then, by Lemma \ref{lem_M_comm}, there is $q$ a type, finitely satisfiable in $M$, such that \begin{equation}\label{mu_tens_q}
    (\mu\otimes q)_{|M}\neq(q\otimes\mu)_{|M}
    \end{equation}
    
    Fix $N$ a saturated extension of $M$. Let $b\models q_{|N}$. Let us build a smooth measure $\nu'_\alpha((z^0_i,z^1_i)_{i<\alpha})$ (with the order $2\alpha$ on the variables, i.e. $z^0_0,z^1_0,z^0_1,\dots$) such that:
    \begin{enumerate}
        \item $(\delta_b\times\nu'_\alpha)_{|N}$ is finitely satisfiable in $M$.
        \item We have $\nu'_{\alpha|N}= F_M(\mu)_{rev|N}^{(2\alpha)}$.
        \item For all $i<\alpha$, we have $(\delta_b\times \nu'_{\alpha|z^0_i})_{|M} \neq (\delta_b\times \nu'_{\alpha|z^1_i})_{|M}$.
    \end{enumerate}

If we build such measures $\nu'_{\alpha}$ extending each other, the process will stop at some $\alpha<(|M|+|T|)^+$. Indeed, if this were not the case, denoting $\nu' := \nu'_{(|M|+|T|)^+}$, we may assume (by the pigeonhole principle) that there are $\phi(x,y)\in L(M)$ and $n\in\N^*$ such that for all $i<(|M|+|T|)^+$, we have $|\nu'(\phi(z^0_i;b))-\nu'(\phi(z^1_i;b))|\geq \frac 1n$. But $\nu'_{|N} = F_M(\mu)^{(2(|M|+|T|)^+)}_{rev|N}$ by construction, hence $\nu'$ is indiscernible over $N$ (by the remarks after Fact \ref{Morley_ind}), so \[\nu'(\phi(z^0_0;b)),\nu'(\phi(z^1_0;b)),\nu'(\phi(z^0_1;b)),\dots\] converges by NIP (Fact\ref{inv_seq_mes}). Contradiction.

So we have shown that the process would stop at some $\alpha<(|M|+|T|)^+$; let us show that we can go one step further anyway. This will give us the desired contradiction. We take $\nu'_\alpha$ such that it is impossible to extend it into $\nu'_{\alpha+1}$ satisfying the desired properties. We write $\nu'$ for $\nu'_\alpha$. If $\lambda(x,y)$ is a measure, recall that we denote by $\lambda_{|x,A}$ the restriction of $\lambda$ to the variable $x$ and the set of parameters $A$. Note that, if $\lambda$ is $M$-invariant, then $F_M(\lambda_{|x}) = F_M(\lambda)_{|x}$.


Let $s_0$ be the unique global extension of $(\delta_b\times\nu')_{|N}$ which is finitely satisfiable in $M$ (it exists by Condition $1$, and is unique by saturation). Denote by $y$ the variable attached to $b$, and by $\overline z = (z_i^0,z_i^1)_{i<\alpha}$ the variable attached to $\nu'$. Let $s = \mu_{x}\otimes s_0$ (we mean that $s= \mu\otimes s_0$, with variable $x$ for the $\mu$ on the left). Then $s$ is $M$-invariant as is a product of $M$-invariant measures. We have:
\begin{equation}\label{s_0_y}
s_{0|y} = q
\end{equation}
Indeed, $s_{0|y,N} = (\delta_b\times\nu')_{|y,N} = \delta_{b|N} = q_{|N}$. But $s_{0|y}$ and $q$ are both $M$-invariant, and $N$ is saturated over $M$, so we get the equality.

Similarly:
\begin{equation}\label{s_0_z}
s_{0|\overline z} = F_M(\mu)_{rev}^{(2\alpha)}
\end{equation}
because $s_{0|\overline z, N} = \nu'_{|N} = F_M(\mu)_{rev|N}^{(2\alpha)}$ by Condition 2. Both sides are $M$-invariant, hence the equality.



\underline{Step $1$}: We prove that $F_M(s)_{|x\overline z} = s_{0|\overline z}\otimes F_M(\mu)$.

Recall that $\mu$ and $F_M(\mu)$ commute. Then, by induction on $n$, $\mu$ and $F_M(\mu)^{(n)}$ commute. Indeed: 
\[\begin{array}{rcl}
\mu(y)\otimes F_M(\mu)^{(n+1)}(x_0,\dots,x_n) & = & \mu(y)\otimes(F_M(\mu)(x_n)\otimes F_M(\mu)^{(n)}(x_0,\dots,x_{n-1})) \\
& = & (\mu(y)\otimes F_M(\mu)(x_n))\otimes F_M(\mu)^{(n)}(x_0,\dots,x_{n-1}) \\
& = & (F_M(\mu)(x_n)\otimes\mu(y))\otimes F_M(\mu)^{(n)}(x_0,\dots,x_{n-1}) \\
& = & F_M(\mu)(x_n)\otimes(\mu(y)\otimes F_M(\mu)^{(n)}(x_0,\dots,x_{n-1})) \\
& = & F_M(\mu)(x_n)\otimes(F_M(\mu)^{(n)}(x_0,\dots,x_{n-1})\otimes \mu(y)) \\
& = & (F_M(\mu)(x_n)\otimes F_M(\mu)^{(n)}(x_0,\dots,x_{n-1}))\otimes \mu(y) \\
& = & F_M(\mu)^{(n+1)}(x_0,\dots,x_n)\otimes\mu(y)
\end{array}\]
So, for all $n<\omega$, we have $\mu\otimes F_M(\mu)^{(n)} = F_M(\mu)^{(n)}\otimes\mu$. But, to check equality between two measures with an infinite number of variables, it is enough to check equality for every finite subset of the variables. So, for any linear order $I$, we have $\mu\otimes F_M(\mu)^{(I)} = F_M(\mu)^{(I)}\otimes\mu$. In particular, by Equation \ref{s_0_z}, we have $s_{|x\overline z} = \mu\otimes s_{|\overline z} = \mu\otimes s_{0|\overline z} = s_{0|\overline z}\otimes \mu$. Now apply $F_M$ to this equality:

\[\begin{array}{rcll}
F_M(s)_{|x\overline z} & = & F_M(s_{|x\overline z}) \\
& = & F_M(s_{0|\overline z}\otimes\mu)\\
& = & s_{0|\overline z}\otimes F_M(\mu)&\text{using Corollary \ref{F_M_comm_prod}}
\end{array}\]

which is the result we anounced.

\underline{Step $2$}: We show that $F_M(s)_{|xy, M}\neq (q\otimes F_M(\mu))_{|M}$.

We have: 
\[\begin{array}{rcll}
s_{|xy, M} & = & (\mu_{x}\otimes s_{0|y})_{|M} &\text{by definition}\\
& = & (\mu\otimes q)_{|M} &\text{by Equation \ref{s_0_y}}\\
&\neq & (q\otimes\mu)_{|M} & \text{by Equation \ref{mu_tens_q}}
\end{array}\] 
Now recall that, for $\lambda$ an $M$-invariant measure, we have $F_M(\lambda)_{|M} = \lambda_{|M}$. So applying $F_M$ yields:

\[\begin{array}{rcll}
F_M(s_{|xy})_{|M} & = & s_{|xy, M}\\
& \neq & (q\otimes\mu)_{|M}\\
& = & F_M(q\otimes\mu)_{|M}\\
& = &(q\otimes F_M(\mu))_{|M}& \text{by Corollary \ref{F_M_comm_prod}}
\end{array}\]

as desired.
 
 \underline{Step $3$}: We finally build our new measure extending $\nu'$.

Note that, as $s_0 = s_{|y\overline z}$ is finitely satisfiable in $M$, we have $F_M(s)_{|y\overline z} = s_0$. In particular, $F_M(s)_{|y\overline z, N} = s_{0|N} = (\delta_b\times\nu')_{|N}$. In other words, $\delta_b\times\nu'$ is a smooth extension of $F_M(s)_{|y\overline z, N}$. So, by Proposition \ref{extension_lem}, there exists $F'$ a smooth extension of $F_M(s)_{|N}$ such that $F'_{|y\overline z} = \delta_b\times\nu'$. Let then be $G'$ a smooth measure such that $G'_{|N} = (F_M(s)\otimes F_M(\mu)_w)_{|N}$ and $G'_{|xy\overline z} = F'$, which exists, again by Proposition \ref{extension_lem}. Then $G'_{|y} = \delta_b$ is a type, so (by the remarks following the definition of a separated amalgam) there exists $\nu''(w,x,\overline z)$ a smooth measure such that $G'=\delta_b\times \nu''$. Set $z_{\alpha}^0:=x$ and $z_{\alpha}^1:=w$ (but keep with $\overline z = (z^0_i,z^1_i)_{i<\alpha}$; this will be a more convenient notation). We clearly have that $\nu''$ extends $\nu'$ (by choice of $G'$ and $F'$). Also, we have that 
\[(\delta_b\times\nu'')_{|N} = G'_{|N} = (F_M(s)\otimes F_M(\mu))_{|N}\]
 is finitely satisfiable in $M$, as a product of measures which are finitely satisfiable in $M$.
 
 Then, note that
 \[\begin{array}{rcll}
 \nu''_{|N} & = & (F_M(s_{|x\overline z})\otimes F_M(\mu))_{|N}&\text{by definition} \\
 & = & (s_{0|\overline z}\otimes F_M(\mu)_x\otimes F_M(\mu)_w)_{|N}&\text{by Step }1\\
 & = & (F_M(\mu)^{(2\alpha)}_{rev}(\overline z) \otimes F_M(\mu)(z^0_\alpha)\otimes F_M(\mu)(z^1_\alpha))_{|N}&\text{by Equation \ref{s_0_z}} \\
 & = & F_M(\mu)^{(2\alpha+2)}_{rev|N} & \text{by the remarks before this theorem}\\
 & = & F_M(\mu)^{(2(\alpha+1))}_{rev|N}
 \end{array}\]
 so condition $2$ is fulfilled. Finally, Condition $3$ follows from Step $2$: indeed,
\[\begin{array}{rcll}
(\delta_b\times \nu''_{|z^1_\alpha})_{|M} & = & G'_{|yw, M}\\
& = & (F_M(s)\otimes F_M(\mu))_{|wy, M}\\
& = & (F_M(s_{|y})\otimes F_M(\mu))_{|M}\\
& = & (F_M(q)\otimes F_M(\mu))_{|M} & \text{by Equation \ref{s_0_y}}\\
& = & (q\otimes F_M(\mu))_{|M} &\text{because }q\text{ is finitely satisfiable in }M\\
& \neq & (\delta_b\times\nu''_{|z^0_\alpha})_{|M} &\text{by Step }2
\end{array}\]which is Condition $3$. But this is a contradiction.
\end{proof}

\newpage

\printbibliography[heading=bibintoc]

\end{document}